\documentstyle[12pt,amssymb]{article}

\begin{document}

\newtheorem{definition}{Definition $\!\!$}[section]
\newtheorem{prop}[definition]{Proposition $\!\!$}
\newtheorem{lem}[definition]{Lemma $\!\!$}
\newtheorem{corollary}[definition]{Corollary $\!\!$}
\newtheorem{theorem}[definition]{Theorem $\!\!$}
\newtheorem{example}[definition]{\it Example $\!\!$}
\newtheorem{remark}[definition]{Remark $\!\!$}

\newcommand{\Label}{\label}
\newcommand{\nc}[2]{\newcommand{#1}{#2}}
\newcommand{\rnc}[2]{\renewcommand{#1}{#2}}

\nc{\Section}{\setcounter{definition}{0}\section}
\renewcommand{\theequation}{\thesection.\arabic{equation}}
\newcounter{c}
\renewcommand{\[}{\setcounter{c}{1}$$}
\newcommand{\etyk}[1]{\vspace{-7.4mm}$$\begin{equation}\Label{#1}
\addtocounter{c}{1}}
\renewcommand{\]}{\ifnum \value{c}=1 $$\else \end{equation}\fi}


\nc{\bpr}{\begin{prop}}
\nc{\bth}{\begin{theorem}}
\nc{\ble}{\begin{lem}}
\nc{\bco}{\begin{corollary}}
\nc{\bre}{\begin{remark}}
\nc{\bex}{\begin{example}}
\nc{\bde}{\begin{definition}}
\nc{\ede}{\end{definition}}
\nc{\epr}{\end{prop}}
\nc{\ethe}{\end{theorem}}
\nc{\ele}{\end{lem}}
\nc{\eco}{\end{corollary}}
\nc{\ere}{\hfill\mbox{$\Diamond$}\end{remark}}
\nc{\eex}{\end{example}}
\nc{\epf}{\hfill\mbox{$\Box$}}
\nc{\ot}{\otimes}
\nc{\bsb}{\begin{Sb}}
\nc{\esb}{\end{Sb}}
\nc{\ct}{\mbox{${\cal T}$}}
\nc{\ctb}{\mbox{${\cal T}\sb B$}}
\nc{\bcd}{\[\begin{CD}}
\nc{\ecd}{\end{CD}\]}
\nc{\ba}{\begin{array}}
\nc{\ea}{\end{array}}
\nc{\bea}{\begin{eqnarray}}
\nc{\eea}{\end{eqnarray}}
\nc{\be}{\begin{enumerate}}
\nc{\ee}{\end{enumerate}}
\nc{\beq}{\begin{equation}}
\nc{\eeq}{\end{equation}}
\nc{\bi}{\begin{itemize}}
\nc{\ei}{\end{itemize}}
\nc{\kr}{\mbox{Ker}}
\nc{\te}{\!\ot\!}
\nc{\pf}{\mbox{$P\!\sb F$}}
\nc{\pn}{\mbox{$P\!\sb\nu$}}
\nc{\bmlp}{\mbox{\boldmath$\left(\right.$}}
\nc{\bmrp}{\mbox{\boldmath$\left.\right)$}}
\rnc{\phi}{\mbox{$\varphi$}}
\nc{\LAblp}{\mbox{\LARGE\boldmath$($}}
\nc{\LAbrp}{\mbox{\LARGE\boldmath$)$}}
\nc{\Lblp}{\mbox{\Large\boldmath$($}}
\nc{\Lbrp}{\mbox{\Large\boldmath$)$}}
\nc{\lblp}{\mbox{\large\boldmath$($}}
\nc{\lbrp}{\mbox{\large\boldmath$)$}}
\nc{\blp}{\mbox{\boldmath$($}}
\nc{\brp}{\mbox{\boldmath$)$}}
\nc{\LAlp}{\mbox{\LARGE $($}}
\nc{\LArp}{\mbox{\LARGE $)$}}
\nc{\Llp}{\mbox{\Large $($}}
\nc{\Lrp}{\mbox{\Large $)$}}
\nc{\llp}{\mbox{\large $($}}
\nc{\lrp}{\mbox{\large $)$}}
\nc{\lbc}{\mbox{\Large\boldmath$,$}}
\nc{\lc}{\mbox{\Large$,$}}
\nc{\Lall}{\mbox{\Large$\forall$}}
\nc{\bc}{\mbox{\boldmath$,$}}
\rnc{\epsilon}{\varepsilon}
\rnc{\ker}{\mbox{\em Ker}}
\nc{\ra}{\rightarrow}
\nc{\ci}{\circ}
\nc{\cc}{\!\ci\!}
\nc{\T}{\mbox{\sf T}}
\nc{\can}{\mbox{\em\sf T}\!\sb R}
\nc{\cnl}{$\mbox{\sf T}\!\sb R$}
\nc{\lra}{\longrightarrow}
\nc{\M}{\mbox{Map}}
\rnc{\to}{\mapsto}
\nc{\imp}{\Rightarrow}
\rnc{\iff}{\Leftrightarrow}
\nc{\bmq}{\cite{bmq}}
\nc{\ob}{\mbox{$\Omega\sp{1}\! (\! B)$}}
\nc{\op}{\mbox{$\Omega\sp{1}\! (\! P)$}}
\nc{\oa}{\mbox{$\Omega\sp{1}\! (\! A)$}}
\nc{\inc}{\mbox{$\,\subseteq\;$}}
\nc{\de}{\mbox{$\Delta$}}
\nc{\spp}{\mbox{${\cal S}{\cal P}(P)$}}
\nc{\dr}{\mbox{$\Delta_{R}$}}
\nc{\dsr}{\mbox{$\Delta_{\cal R}$}}
\nc{\m}{\mbox{m}}
\nc{\0}{\sb{(0)}}
\nc{\1}{\sb{(1)}}
\nc{\2}{\sb{(2)}}
\nc{\3}{\sb{(3)}}
\nc{\4}{\sb{(4)}}
\nc{\5}{\sb{(5)}}
\nc{\6}{\sb{(6)}}
\nc{\7}{\sb{(7)}}
\nc{\hsp}{\hspace*}
\nc{\nin}{\mbox{$n\in\{ 0\}\!\cup\!{\Bbb N}$}}
\nc{\al}{\mbox{$\alpha$}}
\nc{\bet}{\mbox{$\beta$}}
\nc{\ha}{\mbox{$\alpha$}}
\nc{\hb}{\mbox{$\beta$}}
\nc{\hg}{\mbox{$\gamma$}}
\nc{\hd}{\mbox{$\delta$}}
\nc{\he}{\mbox{$\varepsilon$}}
\nc{\hz}{\mbox{$\zeta$}}
\nc{\hs}{\mbox{$\sigma$}}
\nc{\hk}{\mbox{$\kappa$}}
\nc{\hm}{\mbox{$\mu$}}
\nc{\hn}{\mbox{$\nu$}}
\nc{\la}{\mbox{$\lambda$}}
\nc{\hl}{\mbox{$\lambda$}}
\nc{\hG}{\mbox{$\Gamma$}}
\nc{\hD}{\mbox{$\Delta$}}
\nc{\th}{\mbox{$\theta$}}
\nc{\Th}{\mbox{$\Theta$}}
\nc{\ho}{\mbox{$\omega$}}
\nc{\hO}{\mbox{$\Omega$}}
\nc{\hp}{\mbox{$\pi$}}
\nc{\hP}{\mbox{$\Pi$}}
\nc{\bpf}{{\it Proof.~~}}
\nc{\slq}{\mbox{$A(SL\sb q(2))$}}
\nc{\fr}{\mbox{$Fr\llp A(SL(2,\IC))\lrp$}}
\nc{\slc}{\mbox{$A(SL(2,\IC))$}}
\nc{\af}{\mbox{$A(F)$}}
\rnc{\widetilde}{\tilde}
\nc{\qdt}{quantum double-torus}
\nc{\aqdt}{\mbox{$A(DT^2_q)$}}
\nc{\dtq}{\mbox{$DT^2_q$}}
\nc{\uc}{\mbox{$U(2)$}}
\nc{\uq}{\mbox{$U_{q^{-1},q}(2)$}}
\nc{\auq}{\mbox{$A(U_{q^{-1},q}(2))$}}

\def\esl{{\mbox{$E\sb{\frak s\frak l (2,{\Bbb C})}$}}}
\def\esu{{\mbox{$E\sb{\frak s\frak u(2)}$}}}
\def\zf{{\mbox{${\Bbb Z}\sb 4$}}}
\def\zt{{\mbox{$2{\Bbb Z}\sb 2$}}}
\def\ox{{\mbox{$\Omega\sp 1\sb{\frak M}X$}}}
\def\oxh{{\mbox{$\Omega\sp 1\sb{\frak M-hor}X$}}}
\def\oxs{{\mbox{$\Omega\sp 1\sb{\frak M-shor}X$}}}
\def\Fr{\mbox{Fr}}
\def\gal{-Galois extension}
\def\hge{Hopf-Galois extension}
\def\qd{$q$-deformation}
\def\ta{\tilde a}
\def\tb{\tilde b}
\def\tc{\tilde c}
\def\td{\tilde d}
\def\st{\stackrel}

\def\inbar{\,\vrule height1.5ex width.4pt depth0pt}
\def\IC{{\Bbb C}}
\def\IZ{{\Bbb Z}}
\def\IN{{\Bbb N}}
\def\otc{\otimes_{\IC}}
\def\ra{\rightarrow}
\def\ota{\otimes_ A}
\def\otza{\otimes_{ Z(A)}}
\def\otc{\otimes_{\IC}}
\def\h{\rho}
\def\x{\zeta}
\def\th{\theta}
\def\s{\sigma}
\def\t{\tau}


\title{\vspace*{-15mm} \bf QUANTUM DOUBLE-TORUS
\thanks{Dedicated to the memory of Kazimiera Rzeczkowska.}
}

\author{
\normalsize\sc Piotr M.~Hajac\thanks{
On~leave from: Department of Mathematical Methods in Physics, 
Warsaw University, ul.~Ho\.{z}a 74, Warsaw, \mbox{00--682~Poland}.
http://info.fuw.edu.pl/KMMF/ludzie\underline{~~}ang.html~ 
(E-mail: pmh@fuw.edu.pl) Partially supported by the IHES visiting stipend,
NATO postdoctoral fellowship and KBN grant 2 P03A 030 14.
}\\ 
\vspace*{-1mm}
\normalsize Department of Applied Mathematics and Theoretical Physics,\\
\vspace*{-1mm}
\normalsize University of Cambridge, Silver St., Cambridge CB3 9EW, England.\\ 
\vspace*{-1mm}
\normalsize http://www.damtp.cam.ac.uk/user/pmh33\\ 
\vspace*{-1mm}
\normalsize e-mail: pmh33@amtp.cam.ac.uk \\ \ \\
%
\normalsize\sc Tetsuya Masuda\\
\vspace*{-1mm}
\normalsize Institute of Mathematics,\\ 
\vspace*{-1mm}
\normalsize University of Tsukuba, Ibaraki, 305 Japan.\\
\vspace*{-1mm}
\normalsize  e-mail: tetsuya@math.tsukuba.ac.jp\\
\vspace*{-15mm}
}

\date{}
\maketitle
~\\

{
\baselineskip=0pt\small
{\bf Abstract.}
A symmetry extending the $T^2$-symmetry of the noncommutative torus
$T^2_q$ is studied in the category of quantum groups. This extended
symmetry is given by the \qdt\ \dtq\ defined as a compact matrix 
quantum group consisting of the disjoint union of $T^2$ and $T^2_{q^2}$.
The bicross-product  structure of the polynomial Hopf algebra \aqdt\ 
is computed.
The Haar measure and the complete list of unitary irreducible representations
of \dtq\ are determined explicitly.
\\

}

{
\small
{\bf R\'esum\'e.}
Une sym\'etrie qui prolonge la $T^2$-sym\'etrie d'un tore
noncommutatif
$T^2_q$ est \'etudi\'ee dans le contexte des groupes quantiques. Cette
sym\'etrie est donn\'ee par le double-tore quantique \dtq\
d\'efini comme un groupe quantique compact matriciel, 
union disjointe de $T^2$ et $T^2_{q^2}$. La structure de bi-produit
crois\'ee de l'alg\`ebre de Hopf polynomiale \aqdt\ est
calcul\'ee. La mesure d'Haar et la liste compl\`ete des
repr\'esentations irreductibles unitaires de \dtq\ sont determin\'ees
explicitement.
\\ \ \\

}

\centerline{\bf Version fran\c{c}aise abr\'eg\'ee}

\vspace*{2.5mm}

Nous prolongeons la sym\'{e}trie classique $T^2$ du deux-tore noncommutatif
$T^2_q$ \`a la sym\'{e}trie du groupe quantique donn\'e par le double-tore
quantique~\dtq. Contrairement \`a la sym\'{e}trie $T^2$ de $T^2_q$,
la sym\'etrie
\dtq\ d\'{e}couvre la nature quantique des deux-tores noncommutatifs.
Un double-tore quantique est un sous-groupe quantique (non connexe)
de la
$q$-d\'eformation pure tordue de $U(2)$ donn\'ee par la $R$-matrice diagonal
$diag(1,q^{-1},q,1)$. Du point de vue de la th\'eorie de Hopf-Galois, l'anneau
de coordonn\'ees de \dtq\ est un bi-produit crois\'e de l'alg\`ebre d'Hopf
polynomiale de $Z_2$ et de celle de $T^2$ avec un cocycle non-trivial
et une coaction qui d\'eterminent les structures d'alg\`ebre et de 
cog\'ebre, respectivement. (Le groupe quantique \dtq\ peut \^{e}tre vu comme une 
quantification du
produit semi-direct du groupe~$T^2\rtimes Z_2$.) La completion de
l'anneau de coordonn\'ees de \dtq\ nous donne une $C^*$-alg\`ebre qui est
isomorphe \`a la somme directe de l'alg\`ebre $C(T^2)$ des fonctions continues
sur le deux-tore classique et la $C^*$-alg\'ebre $C(T^2_{q^2})$ de
deux-tore noncommutatif. \dtq\ est un groupe quantique compact matriciel
dans le sens de S.~L.~Woronowicz. Son dual unitaire (collection des toutes les
repr\'esentations irreductibles unitaires) est donn\'e par la formule:
$
\widehat{\dtq}={\Bbb Z}\coprod {\Bbb Z}\coprod ({\Bbb Z}\times {\Bbb N}_{+}).
$

\section*{Introduction}

The aim of this note is to construct and analyse a $q$-deformation of the
semi-direct product group $T^2\rtimes Z_2$ that extends the $T^2$-symmetry
of the noncommutative torus $T^2_q$. This semi-direct product is a subgroup
of $U(2)$ consisting of diagonal and off-diagonal matrices with
unitary entries.
Since it is a disjoint union of two tori, we call it a double torus, and 
denote~$DT^2$.

To carry out a $q$-deformation of $DT^2$ it is convenient to start with
the well-known two-parameter quantisation $GL_{p,q}(2)$, $p,q\in{\Bbb C}^\times$,
 of $GL(2)$ given by the $R$-matrix
\beq\Label{r}
\pmatrix{
1 & 0 & 0 & 0 \cr
0 & q^{-1} & 0 & 0 \cr
0 & 1-(pq)^{-1} & p^{-1} & 0 \cr
0 & 0 & 0 & 1 \cr
}.
\eeq
For $p=\bar{q}$ we can consider the involutive structure given on generators
by
$
a^*=S(a),~ b^*=S(c),~ c^*=S(b),~ d^*=S(d)
$,
where $S$ is the antipode. The introduction of such an involutive structure
makes the quantum determinant $D:=ad-\bar{q}bc=da-q^{-1}bc$ a unitary element,
and defines the coordinate ring $A(U_{\bar{q},q}(2))$ of the compact matrix
quantum group $U_{\bar{q},q}(2)$. (See \cite{mhr} for details.)

For $q\in{\Bbb R}^\times$, the quantum determinant $D$ becomes central, 
we can set it to 1, and decompose 
$U_{q,q}(2)$ as follows: $U_{q,q}(2)=:U_q(2)=T^1\times SU_q(2)$.

The other characteristic case is when $q\in U(1)$. Then the $R$-matrix~(\ref{r})
becomes diagonal, and the commutation relations between generators of \auq\
are
$
ab=q^{-1}ba,~ cd=q^{-1}dc,~ ac=qca,~ bd=qdb,~ ad=da,~ bc=q^2cb\, 
$.
This \qd\ of $U(2)$ is precisely what we use to quantise the double-torus
subgroup of~$U(2)$.

Although it is not technically needed, we assume the coefficient field to be
$\IC$ and introduce the involutive structure from the very beginning for the
sake of brevity and clarity. For the coaction and coproduct we use the Sweedler
notation with the suppressed summation sign.

\section{Quantum Double-Torus as a Symmetry of \boldmath $T^2_q$}

One can represent the 2-torus $T^2$ as the quotient ${\Bbb R}^2/{\Bbb Z}^2$
or as a subset of ${\Bbb C}^2$. The latter point of view leads to
the question of which subgroups of $GL(2,{\Bbb C})$ preserve $T^2$. A natural
subgroup of $GL(2,{\Bbb C})$ that preserves $T^2$ is the double
torus $DT^2$ described in the Introduction. The double torus is
a subgroup of $U(2)$, which can be viewed as a field of $T^3$
 over the internal points of $[0,1]$ bounded by $T^2$ at each of 
the endpoints.  
Although this geometrical picture disregards the Lie group structure
of $U(2)$, it helps us to describe the support of the Haar functional on the Hopf
algebra $A(U(2))$ of the polynomial functions on $U(2)$.
This support is included in the space of polynomials on the aforementioned
interval $[0,1]$, and the Haar functional is given by integration on~$[0,1]$.
Moreover, this picture of $U(2)$ allows us to locate $DT^2$ inside $U(2)$ as
the two $T^2$ placed at the endpoints of~$[0,1]$. What we study here is a
$q$-deformed version of this setting.

The standard $q$-deformation of $U(2)$ turns the interval $[0,1]$ into the
quantum interval $[0,1]_q$ consisting of powers of~$q$. Missing is the
 endpoint corresponding to the off-diagonal torus of $DT^2$. The 
off-diagonal part of the double-torus subgroup of \uc\ is destroyed by
this quantisation. On the other hand, the pure-twist $q$-deformation of
\uc, with $q\in U(1)$, leaves the interval $[0,1]$ unchanged. It deforms
\uc\ into \uq, which can be understood as a field of {\em
noncommutative} 3-tori over the open interval $(0,1)$ bounded by $T^2_{q^2}$ 
at 0 and $T^2$ at 1 (see~\cite[Section~6]{r-m}, cf.~\cite{m-k91}). 
The union of $T^2_{q^2}$ and $T^2$ is the {\em quantum double-torus}
subgroup of \uq. More precisely, we have:
\bde\Label{qdtdef}
Let $I$ be the Hopf ideal of \auq\ generated by $\{ab,ac,cd,bd\}$. We call
the quotient Hopf algebra $\aqdt:=\auq/I$ the coordinate ring of the {\em
\qdt}~\dtq.
\ede

\aqdt\ is a $*$-Hopf algebra. The Diamond Lemma \cite[Theorem~1.2]{b} enables us
to claim the following:
\ble\Label{diamond}
Put $z=D^{-1}ad$. The set of monomials:
\begin{eqnarray*}
&&
D^kz^la^mc^n,~ m>0,~ n>0;~~ 
D^kz^la^mb^{-n},~ m>0,~ n<0;~~ 
\\ &&
D^kz^ld^{-m}c^n,~ m<0,~ n>0;~~
D^kz^ld^{-m}b^{-n},~ m<0,~ n<0;~~
\end{eqnarray*}\vspace*{-9mm}\begin{eqnarray*}
&&
D^kz^lc^{n},\; n>0;~~~ 
D^kz^lb^{-n},\; n<0;~~~
\\ &&
D^kz^la^m,\; m>0;~~~ 
D^kz^ld^{-m},\; m<0;~~~
\end{eqnarray*}\vspace*{-8mm}\beq\Label{b1} 
D^kz^l;~~~k,l,m,n\in\IZ, 
\eeq
is a $\IC$-basis of \auq. Denote the generators of the quotient algebra
\aqdt\ by the same letters as generators of~\auq. The set of polynomials
\beq\Label{b2}
D^ma^n,~~ D^md^n,~~ D^mz,~~ D^mb^n,~~ D^mc^n,~~ D^m(1-z);~~~ m,n\in\IZ,~ n>0;~~~
\eeq
is a $\IC$-basis of \aqdt.
\ele

\bre\Label{discon}\em
It is tempting to define a {\em disconnected quantum space} by requiring its
function algebra to contain a non-trivial {\em central} idempotent. 
Since $z$ is central and becomes an
idempotent after passing to the quotient algebra~\aqdt, the 
\qdt\ is disconnected
according to such a definition.
\ere

The main point of this section is that \dtq\ acts the same way on $T^2_q$
as $DT^2$ on~$T^2$. Thus we obtain a {\em quantum symmetry} of $T^2_q$ that
extends the well-studied classical $T^2$-symmetry of~$T^2_q$. More precisely,
we have:
\bpr\Label{diagram} 
Let $\rho_c$ and $\rho_q$ be the appropriate left coactions induced by the
standard left coaction given on generators by 
{\scriptsize
$\pmatrix{x \cr y \cr}\to
\pmatrix{
a\ot 1 & b\ot 1 \cr
c\ot 1 & d\ot 1 \cr
}\pmatrix{1\ot x \cr 1\ot y \cr
}$}.
Let $\pi$ be the 
canonical projection setting $z=1$. Then the following diagram is commutative:
\[
\begin{array}{ll}
A(T^2_q) & \mathop{-\hspace{-6pt}\longrightarrow}\limits^{\rho_c} 
~~A(T^2)\ot A(T^2_q)  \\ \ \\
~~^{id} \Big\uparrow  & ~~~~~~~~~~~~~~^{\pi\ot id}\Big\uparrow\  \\ \ \\
A(T^2_q) & \mathop{-\hspace{-6pt}\longrightarrow}\limits^{\rho_q}
~~\aqdt\ot A(T^2_q)\, .\\
\end{array}
\]
\epr
 A standard reasoning allows us to claim 
the same at the $C^*$-algebraic level.

The classical $T^2$-symmetry of $T^2_q$ is insensitive to
the noncommutative nature of $T^2_q$. In contrast, 
the \dtq-quantum symmetry of $T^2_q$ detects
the square of the deformation parameter $q$, i.e., we have $bc=q^2cb$ in \aqdt.

\section{Bicross-product Structure of \boldmath\aqdt}

The concept of short exact sequences of groups can be extended to the category
of quantum groups, where it is described by (strictly) exact sequences of
Hopf algebras (see~\cite{pw,s2,ad}).  
Consider the sequence of Hopf algebras
$
A(Z_2)\st{i}{\ra}\aqdt\st{\pi}{\ra}A(T^2)
$,
where $i$ is the Hopf algebra injection sending $f\in A(Z_2)$, $f(0):=1$, 
$f(1):=0$, to $z$, and $\pi$ is obtained from the canonical surjection
setting $z=1$. It can be directly verified that this
 is a strictly exact sequence of Hopf algebras in the sense
of~\cite[p.3338]{s2}. Thus we can view \dtq\ as a quantum-group extension
of $Z_2$ by~$T^2$. In the language of the Hopf-Galois theory, we can describe
\aqdt\ as a cleft $A(T^2)$\gal\ of $A(Z_2)$, or a crossed product algebra 
of $A(Z_2)$ by~$A(T^2)$ (e.g., see~\cite{s3}).

Recall that a cleaving map of a cleft $H$\gal\ $B\inc P$ is a convolution
invertible colinear map from a Hopf algebra $H$ to an $H$-comodule 
algebra~$P$. In the case of \aqdt\ we can construct a cleaving map
$j:A(T^2)\ra\aqdt$ in the following way:
\beq\Label{j}
j(u^kv^l):=\left\{
\begin{array}{ll}
D^la^{k-l}+(-1)^lq^{l^2}D^lc^{k-l} & \mbox{for $k>l$}\\
D^kz+(-D)^{-k}(1-z) & \mbox{for $k=l$}\\
(-1)^kq^{-k^2}D^kb^{l-k}+D^kd^{l-k} & \mbox{for $k<l$}. 
\end{array}
\right.
\eeq
Here $u,v$ are generators of $A(T^2)$, and $k,l\in\IZ$. Observe that for
$q\neq 1$ the map $j$ is a $*$-homomorphism but not an algebra map. It is
an algebra homomorphism only for~$q=1$. Then it can be obtained as the
pull-back of the appropriate trivialisation map for the trivial principal
bundle $DT^2(Z_2,T^2)$. We can think of \dtq\ as a quantum principal bundle
with the classical base space $Z_2$ and classical structure group $T^2$,
but whose total space consists of two non-equivalent fibres: $T^2$ 
and~$T^2_{q^2}$.

It is known that cleft \hge s and crossed products are equivalent notions
(see Theorem~1.12 in~\cite{s3}, \cite{bm,dt1}). The crossed-product structure
of \aqdt\ is given by a {\em non-trivial} cocycle 
$\hs:A(T^2)\ot A(T^2)\ra A(Z_2)$ and the trivial action of $A(T^2)$ on~$A(Z_2)$.
One can compute $\hs$ according to the formula (cf.~\cite[Lemma~1.9]{dhs97}):
$
\hs(h\ot g)=j(h\1)j(g\1)j^{-1}(h\2 g\2)
$.
Explicitly, we obtain: 
\bpr\Label{cocycle}
The cocycle $\hs$ associated to the cleaving map~(\ref{j}) is given by
the formulas:
\[
\hs=\hs_cz+\hs_q(1-z)\, ,~~~ \hs_c=\he\ot\he~~~
\mbox{($\he$ is the counit map),}
\]{
\[
\hs_q(u^kv^l\ot u^mv^n)=\left\{
\begin{array}{ll}
q^{-2kn} & \mbox{for $k>l$, $m>n$}\\
q^{-m(2k+m)} & \mbox{for $k>l$, $m=n$}\\
q^{-2n(k+m)} & \mbox{for $k>l$, $m<n$, $k+m>l+n$}\\
q^{(k+m)(2l-k-m)} & \mbox{for $k>l$, $m<n$, $k+m=l+n$}\\
q^{2l(k+m)} & \mbox{for $k>l$, $m<n$, $k+m<l+n$}\\
q^{-k(k+2n)} & \mbox{for $k=l$, $m>n$}\\
1 & \mbox{for $k=l$, $m=n$}\\
q^{k(k+2m)} & \mbox{for $k=l$, $m<n$}\\
q^{-2k(l+n)} & \mbox{for $k<l$, $m>n$, $k+m>l+n$}\\
q^{(l+n)(l+n-2k)} & \mbox{for $k<l$, $m>n$, $k+m=l+n$}\\
q^{2m(l+n)} & \mbox{for $k<l$, $m>n$, $k+m<l+n$}\\
q^{m(m+2l)} & \mbox{for $k<l$, $m=n$}\\
q^{2lm} & \mbox{for $k<l$, $m<n$}~. 
\end{array}
\right.
\]}
\epr
The cocycle $\hs$ determines an algebra structure on $A(Z_2)\ot A(T^2)$
by the following rule (see~\cite[p.693]{bcm}):
$
(x\ot h)(y\ot g)=xy\hs(h\1\ot g\1)\ot h\2 g\2\,
$.
The tensor product
$A(Z_2)\ot A(T^2)$ with the above algebra structure is isomorphic as an algebra
to \aqdt.

One can determine the crossed coproduct structure of~\aqdt\ in a
similar fashion. First,
we can transform a cleaving map $j$ into a cocleaving map $\ell$ given by the
formula $\ell(p)=p\0 j^{-1}(p\1)$. (The inverse transform is given 
by~\cite[3.2.13(2)]{ad}.) The mapping $\ell$ associated to (\ref{j}) can be 
 computed directly to be:
\[
\ell(D^ma^n)=z,~~~ 
\ell(D^md^n)=z,~~~
\ell(D^mz)=z,~~~
\ell(D^m(1-z))=(-1)^m(1-z),
\]\beq\Label{l}
\ell(D^mb^n)=(-1)^mq^{m^2}(1-z),~~~  
\ell(D^mc^n)=(-1)^mq^{-m^2}(1-z).~ 
\eeq
Here $m,n\in\IZ,~ n>0$. 
It can be verified that $\ell$ is a $*$-homomorphism, and that it coincides
with its convolution inverse.

The crossed coproduct structure of \aqdt\ is given by a {\em non-trivial}
coaction $\hl:A(T^2)\ra A(T^2)\ot A(Z_2)$, and the trivial cococycle.
We calculate $\hl$ from $\ell$ according to the 
formula~\cite[Proposition~3.2.9]{ad}:
$
\hl(\pi(p))=\pi(p\2)\ot\ell^{-1}(p\1)\ell(p\3)\, 
$.
Hence we obtain:
\bpr\Label{coact}
The coaction $\hl$ associated to the cocleaving map 
(\ref{l}) is a $*$-algebra homomorphism defined by
$
\hl(u^mv^n)=u^mv^n\ot z+u^nv^m\ot(1-z).
$
\epr
For $q=1$, the coaction $\hl$ is the pull-back of the $Z_2$ action on $T^2$
defining $DT^2$ as the semi-direct product. The fact that $\hl$ remains
unchanged under the \qd\ reflects the property that matrix quantum groups,
in particular \dtq, have the classical (matrix) coproduct. The coaction $\hl$ 
determines a coalgebra structure on $A(Z_2)\ot A(T^2)$ by the following rule
(see~\cite[(6)]{m-s97}):
$
\hD(x\ot h)=(x\1\ot {h\1}^{(0)})\ot(x\2 {h\1}^{(1)}\ot h\2)
$,
where $\hl(h):=h^{(0)}\ot h^{(1)}$. The tensor product
$A(Z_2)\ot A(T^2)$ with the above 
coalgebra structure is isomorphic as a coalgebra to \aqdt.

The crossed-product and crossed-coproduct structures described above together
form a bicross-product (see~\cite[p.244]{m-s97} for this type of bicross-products):
\bpr\Label{bic}
\aqdt\ is isomorphic as a Hopf algebra to the bicross-product
$A(Z_2)^{\lambda}\!\#_\sigma A(T^2)$
defined by cocycle \hs\ and coaction \hl\ of propositions \ref{cocycle} and
\ref{coact}, respectively.
\epr
The same kind of bicross-product structure (i.e., with the trivial action and
cococycle) is studied for the Borel quantum subgroup of $SL_{e^{2\pi i/3}}(2)$ 
in~\cite{dhs97}, and affine Hopf algebra $U_q(\hat{sl_2})$ in~\cite{m-s97}.

To end this section, let us mention that for any $q\in U(1)$ satisfying
$(-q)^{n^2}=1,\,n\in\IZ,\, n>0$,
 we can construct a finite dimensional Hopf algebra, by
dividing \aqdt\ by the Hopf ideal generated by:
$
a^n-z,~ b^n-(1-z),~ c^n-(1-z),~ d^n-z,~ D^n-1 
$.

\section{Representation Theory of \boldmath\dtq}

Since the \qdt\ \dtq\ is a compact quantum group, all of the unitary irreducible
representations are finite dimensional. They are obtained by the irreducible
decompositions of the (multiple) tensor products of some basic representations.
The obvious representation is the unitary group-like element given by the 
(non-central) quantum determinant~$D$. Other basic representations are given
by the (central) unitary group-like element $2z-1$, and the natural vector 
representation {\scriptsize
$\pmatrix{
a & c \cr
b & d \cr
}$}.
There exist the following three families of unitary irreducible representations:
\beq\Label{list}
\chi_m:=D^m,~~ \chi_m^z:=D^m(2z-1),~~ w_{m,n}:=\pmatrix{
D^ma^n & D^mb^n \cr
D^mc^n & D^md^n \cr
},~~~ m,n\in\IZ,~ n>0.
\eeq

The uniquely determined Haar functional $h$ on \aqdt\ is positive and faithful.
It vanishes on all elements of the $\IC$-basis (\ref{b2}) of \aqdt\ except
for the central idempotents $z$ and $(1-z)$, corresponding to the classical
and quantum component of \aqdt, respectively. (It annihilates
all polynomials of $\chi_m$, $\chi_m^z$ and $w_{m,n}$.) 
Let $\{e^{c}_{m,n},e^{q}_{m,n}\}_{m,n\in{\Bbb Z}}$ be an orthonormal  basis of
Hilbert space $\cal H$. The GNS construction
yields the following faithful $*$-representation $\pi$ of \aqdt\ on $\cal H$:
\[
\pi(a)e^{c}_{m,n}=
\left\{
\begin{array}{ll}
e^{c}_{m,n+1} & \mbox{for $n\geq 0$}\\
e^{c}_{m+1,n+1} & \mbox{for $n<0$}, 
\end{array}
\right.
~~~
\pi(b)e^{q}_{m,n}=
\left\{
\begin{array}{ll}
-q^{2n-1}e^{q}_{m+1,n-1} & \mbox{for $n>0$}\\
q^{2m}e^{q}_{m,n-1} & \mbox{for $n\leq 0$}, 
\end{array}
\right.
\]\[
\pi(c)e^{q}_{m,n}=
\left\{
\begin{array}{ll}
e^{q}_{m,n+1} & \mbox{for $n\geq 0$}\\
-q^{-2m-1}e^{q}_{m+1,n+1} & \mbox{for $n<0$}, 
\end{array}
\right.
~~~
\pi(d)e^{c}_{m,n}=
\left\{
\begin{array}{ll}
e^{c}_{m+1,n-1} & \mbox{for $n>0$}\\
e^{c}_{m,n-1} & \mbox{for $n\leq 0$}, 
\end{array}
\right.
\]
$\pi(a)e^{q}_{m,n}=0,~ \pi(b)e^{c}_{m,n}=0,~ \pi(c)e^{c}_{m,n}=0,~ 
\pi(d)e^{q}_{m,n}=0$.
The operator norm completion $C(\dtq)$ of \aqdt\ is a $C^*$-algebra isomorphic
to the direct sum of the algebra $C(T^2)$ of continuous functions on $T^2$
and the noncommutative $C^*$-algebra $C(T^2_{q^2})$, $q\in U(1)$,
 of the quantum two-torus. This means, in particular, that the $C^*$-algebra
$C(\dtq)$ is not of type I for a generic $q\in U(1)$.

Standard reasoning allows us to conclude that the comultiplication extends
from \aqdt\ to $C(\dtq)$, and we have:
\bpr\Label{cmqg}
\dtq\ is a compact matrix quantum group in the sense of
\cite[Definition~1.1]{w-s}.
\epr
\bth\Label{the}
The unitary dual of the \qdt\ is given by 
$
\widehat{\dtq}={\Bbb Z}\coprod {\Bbb Z}\coprod ({\Bbb Z}\times {\Bbb N}_{+})
$.
The complete list of all unitary irreducible representations of \dtq\
is given by the three families of representations~(\ref{list}).
\ethe


\end{document}